\documentclass[12pt]{jsg}

\input xy
\xyoption{all}
\usepackage{latexsym}
\usepackage{amsmath}
\usepackage{amssymb}
\usepackage{amscd}
\usepackage{epsfig}

\newtheorem{theorem}{Theorem}[section]

\newtheorem{corollary}[theorem]{Corollary}
\newtheorem{lemma}[theorem]{Lemma}
\newtheorem{prop}[theorem]{Proposition}
\newtheorem{definition}[theorem]{Definition}
\newtheorem{claim}[theorem]{Claim}
\newtheorem{remit}[theorem]{Remark}
\newenvironment{rem}{\begin{remit}\rm}{\end{remit}}

\numberwithin{equation}{section}

\def\C{\ensuremath{\mathbb{C}}}
\def\iso{\cong}

\def\t{\mathfrak{t}}
\newcommand{\res}{{\rm Res}}
\newcommand{\liet}{{\mathfrak{t}}}
\newcommand{\liets}{\liet^*}
\newcommand{\Jac}{\bigtriangleup}

\begin{document}

\title{Distinguishing the Chambers of the Moment Polytope}

\author{R.\ F.\ Goldin}
\address[R.\ F.\ Goldin]{
Department of Mathematical Sciences\\ 
 George Mason University MS\# 3F2\\ Fairfax, VA 22030}
\email{rgoldin@math.gmu.edu} 

\author{T.\ S.\ Holm}
\address[T.\ S.\ Holm]{Department of Mathematics\\ 
University of California\\ 
Berkeley, CA 94720-3840}
\email{tsh@math.berkeley.edu}

\author{L.\ C.\ Jeffrey}
\address[L.\ C.\ Jeffrey]{Department of Mathematics\\ University of Toronto\\
Toronto, Ontario, Canada M5S 3G3}
\email{jeffrey@math.toronto.edu}

\date{\today}
\thanks{R.G. was supported in part by NSF-DMS, T.H. was supported by an NSF postdoctoral fellowship, and L.J. was supported in part by NSERC}

\begin{abstract}
Let $M$ be a compact manifold with a Hamiltonian $T$ action and moment
map $\Phi$. The restriction map in rational equivariant cohomology from $M$ to
a level set $\Phi^{-1}(p)$ is a surjection, and we denote the kernel
by $I_p$. When $T$ has isolated fixed points, we show that $I_p$
distinguishes the chambers of the moment polytope for $M$. In
particular, counting the number of distinct ideals $I_p$ as $p$ varies
over different chambers is equivalent to counting the number of
chambers.
\end{abstract}

\maketitle

\section{Introduction}

Let $(M,\omega)$ be a compact symplectic manifold with a Hamiltonian
torus action by $T$. Then $M$ has a moment map 
$$
\Phi:M\longrightarrow\t^*
$$ 
from $M$ to the dual of the Lie algebra of $T$ which has the property
that for every $\xi\in\t$, $d\langle \Phi(p),\xi\rangle= \omega(X_\xi,
\cdot)$, where $X_\xi$ is the vector field on $M$ generated by
$\xi$. The torus action preserves level sets of $\Phi$. If $p$ is a
regular value of the moment map, then $T$ acts locally freely on
$\Phi^{-1}(p)$ and the {\em symplectic reduced space}
$M_p:=\Phi^{-1}(p)/T$ has at worst orbifold singularities. In
particular, the equivariant cohomology of the level set is
(rationally) isomorphic to the ordinary 
cohomology of the
quotient $\Phi^{-1}(p)/T$.
The {\em Kirwan map} 
$$
\kappa_p:H^*_T(M)\longrightarrow
H^*(M_p)
$$ is the restriction to the level set in equivariant cohomology
followed by this isomorphism. Kirwan \cite{Ki:Book} showed that $\kappa_p$ is
a surjection for all regular values $p$. Let 
$$
I_p:=\ker\kappa_p
$$
be the kernel of this map. The main aim of this article is to observe
that this ideal distinguishes connected components of the set of
regular values of $\Phi$, when the fixed point set consists of
isolated points. 

In general, this ideal may be hard to compute. S. Tolman and
J. Weitsman described a set of   generators of $I_p$ in terms of their
restriction  properties to the fixed point set of the $T$ action on
$M$ \cite{TW:abelianreduction}. The first author (of this article)
subsequently described these classes in terms of the moment map
polytope $\Phi(M)$\cite{Go:Effective}. Preceding these techniques, the
third author described the ideal indirectly by using cohomology
pairings and the residue formula \cite{JK:LocalizationNonabelian}, a
technique we exploit further here. 

Let $\Delta:=\Phi(M)$ be the image of the moment map. For $M$ compact,
$\Delta$ is a convex polytope which we call the {\em moment polytope}
\cite{At:convexity}, \cite{GS:convexity}. The connected components of
the set of regular values of $\Phi$ form {\em chambers} of the
polytope. These chambers are bounded by the (non-regular) values
$\bigcup\Phi(M^H)$, where the union is over all one-dimensional
subtori $H\subset T$ and $M^H$ denotes the fixed point set of $H$. For
any particular $H\subset T$, we call the connected components of
$\Phi(M^H)$ {\em walls} of $\Delta$, and note that a wall will be
codimension-1 in $\t^*$ if $H$ is one-dimensional and the quotient
$T/H$ acts effectively on $M^H$. While there have been numerous
attempts to describe the geometry as a regular value $p$ crosses a
codimension-1 wall (most notably \cite{GK:residue},
\cite{Ma:wallcrossing}, \cite{GS:birtationalequivalence}), there has
so far been little evidence that global observations about $\Delta$
can be made using the techniques of equivariant cohomology. In
\cite{GLS:symplecticfibrations} the authors are able to use explicit calculations to find the number of chambers in the image of certain
degenerate $SU(4)$ coadjoint orbits, and E. Rassart
\cite{Ra:numberchambers} has since extended these methods to count the
number for any coadjoint orbit of $SU(4)$. The result we present below
suggest that an algebraic approach may have merit in answering these
questions. 
 
\begin{theorem}\label{th:main} 
Let $M$ be a compact connected symplectic manifold with a Hamiltonian
$T$ action.  Suppose the action has isolated fixed points, and denote
the moment map by $\Phi$. Let $p$ and $q$ be two regular values of
$\Phi$. Then $p$ and $q$ in are in distinct chambers of the moment
polytope if and only if $I_p\not=I_q$.  
\end{theorem}

One of the difficulties in proving such a theorem is that, for $p$ and
$q$ in distinct chambers, it is not sufficient to find different
generating sets for the two ideals. These two sets could generate the
same ideal. For this reason, our task is to find specific classes
which are in one ideal but not in the other.

\begin{rem}
The quotient of $H^*_T(M)$ by each of these ideals may
result in isomorphic rings. For example, the reduced space of any
generic $SU(3)$ coadjoint orbit is a 2-sphere for any regular
value. Thus the reduced spaces have isomorphic cohomology even though
the points of reduction may not be in the same chamber. 
\end{rem}

Before continuing, we would like to thank Etienne Rassart and Victor Guillemin for several insightful conversations with each of the authors. We also thank
the referee for useful remarks on the original manuscript.

\section{Background}

The kernel of the Kirwan map can be described in two equivalent ways:
in terms of cohomology classes and in terms of cohomology
pairings. The Tolman-Weitsman theorem (\cite{TW:abelianreduction}, Theorem 3)  describes  a generating set of
the classes in the kernel of the Kirwan map, while the Jeffrey-Kirwan
theorem describes how to evaluate pairings (or integrals) on the
reduced space. Because the reduced space is finite-dimensional, and in
particular the Poincar\'e pairing is non-degenerate, the kernel $I_p$
of the Kirwan map $\kappa_p$ consists of those classes $\beta$ for
which $\int_{M_p}\kappa_p(\alpha\beta)=0$ for all $\alpha$. Our
interest is to show that certain classes have non-zero integral and
are therefore {\em not} in the kernel of the Kirwan map. 

Our main tool for describing any class is its restriction to the fixed
point set $M^T$. 

\begin{theorem}[Kirwan]\label{th:injectivity} 
Let $M$ be a compact symplectic manifold with Hamiltonian $T$ action
and fixed point set $M^T$. Let $M^T_{cc}$ indicate the connected
components of the fixed point set. The restriction to the fixed point
set 
\begin{equation}\label{eq:injectivity}
\iota^*:H^*_T(M)\longrightarrow H^*_T(M^T) = \oplus_{F\in M^T_{cc}} H^*_T(F)
\end{equation}
is an injection.
\end{theorem}

In particular, a class is determined by its behavior on the fixed
point set. For this reason, we refer to the set of connected
components of the fixed point set where the class is nonzero as the
{\em support} of a class. We denote the restriction of $\alpha$ to
$F\in M^T_{cc}$ by $\iota^*_F(\alpha)$ or $\alpha|_F$. 

\begin{definition}\label{def:support}
For $\alpha\in H_T^*(M)$, the {\em support} of $\alpha$ is 
$$
supp\ \alpha = \{F\in M^T_{cc}:\ \iota^*_F(\alpha)\neq 0\},
$$ 
the set of fixed point components $F$ on which $\alpha$ restricts
to a non-zero class. 
\end{definition}

\subsection{The Jeffrey-Kirwan residue formula} 
The residue formula expresses cohomology pairings on reduced spaces in
terms of a multi-dimensional residue of certain rational holomorphic
functions on $\mathfrak{t}$ \cite{JK:LocalizationNonabelian}. The
functions are obtained by restriction of equivariant classes to fixed
points, and dividing by equivariant Euler classes. While the theorem
applies to compact $M$ reduced by any compact group $K$, we need (and
present) only the Abelian case below.

The computation of terms in the residue formula depends on a cone
$\Lambda$ in $\mathfrak{t}$, even though the result of the formula is
independent of this choice. Let $\gamma_1,\dots,\gamma_k$ be the set
of all weights that occur by the $T$ action at any of the fixed point
components. Choose some $\xi\in\t$ such that $\gamma_i(\xi)\neq 0$ for
all $i$. Let $\beta_i = \gamma_i$ if $\gamma_i(\xi)>0$ and
$\beta_i=-\gamma_i$ if $\gamma_i(\xi)<0$. Thus $\beta_i(\xi)>0$ for
all $i$. The cone $\Lambda$ is the set of all vectors in $\t$ which
behave like $\xi$:
$$
\Lambda = \{X\in\t:\ \beta_i(X)>0,\ \mbox{for all}\ i\}.
$$

\begin{theorem}[Jeffrey-Kirwan] 
Let $(M,\omega)$ be a compact symplectic
manifold with a Hamiltonian $T$ action and moment map $\Phi$, where
$T$ is a compact torus. Denote by $M^T_{cc}$ the connected components of the fixed point set of $T$ on $M$. Let $p$ be a regular value of $\Phi$  and
$\omega_p$ the Marsden-Weinstein reduced symplectic form on
$M_p$. Then for $\beta\in H^*_T(M)$ and $\kappa_p:H_T^*(M)\rightarrow
H^*(M_p)$ we have 
$$
\int_{M_p}\kappa_p(\beta)e^{\omega_p} = c\cdot \res^\Lambda \left(\sum_{F\in
M^T_{cc}} e^{i(\Phi(F)-p)(X)}\int_F\frac{\iota^*_F (\beta(X)
e^{\omega})}{e_F(X)}[dX]\right) 
$$
where $c$ is a non-zero constant, $X$ is a variable in
$\mathfrak{t}\otimes \C$, and $e_F(X)$ is the equivariant Euler class
of the normal bundle to $F$ in $M$. The multi-dimensional residue
$\res^\Lambda$ is defined below. 
\end{theorem}

It is often more useful to refer to the dual cone $\Lambda^*\subset
\t^*$, the convex cone formed by the positive span of the vectors  
$\beta_1,\dots, \beta_k$, or the set of vectors $\beta$ 
which can be written $\beta=\sum_{i=1}^k s_i \beta_i$, with $s_i\geq 0$, 
perhaps not uniquely.  Assume $\dim T=l$, let 
$J=(j_1,\dots,j_l)$ be a multi-index, 
and $X^J = X_1^{j_1}\cdots X_l^{j_l}$. 
The operator ${\rm Res}^\Lambda$ 
(defined in \cite{JK:LocalizationNonabelian}) 
is defined by linearity and the following properties 
(\cite{JK:LocalizationQuantization}, Proposition 3.2):
\begin{enumerate}
\item Let $\alpha_1,\dots,\alpha_v\in \Lambda^*$ be vectors in the dual cone.
Suppose that
 $\lambda$ is not in any cone of 
dimension $l-1$ or less spanned by a subset of the $\{\alpha_i\}$.
 Then  $$
{\rm
Res}^\Lambda\left(\frac{X^Je^{i\lambda(X)}[dX]}{\prod_{i=1}^{v}\alpha_i(X)}\right)=0
$$ 
unless all of the following properties are satisfied:
\begin{enumerate}
\item $\{\alpha_i\}_{i=1}^v$  span $\t^*$ as a vector space,
\item $v-(j_1+\dots +j_l)\geq l$,
\item $\lambda\in \langle \alpha_1,\dots, \alpha_v\rangle^+,$ the
positive span of the vectors $\{\alpha_i\}.$ 
\end{enumerate}

\item If properties (1)(a)-(c) above are satisfied, then 
$$
{\rm
Res}^\Lambda\left(\frac{X^Je^{i\lambda(X)}[dX]}
{\prod_{i=1}^{v}\alpha_i(X)}\right)
= \sum_{m\geq 0}\lim_{s\to 0^+}  
{\rm Res}^\Lambda\left(\frac{X^J (i\lambda(X))^m 
e^{is\lambda(X)}[dX]}{m!\prod_{i=1}^{v}\alpha_i(X)}\right),$$
and all but one term in this sum are $0$ (the non-vanishing term being
that with $m = v-( j_1+\dots +j_l)     -l$).
\item \label{it:nonzeroresidue} 
The residue is not identically $0$. If properties (1) $(a)-(c)$ are
satisfied with $\alpha_{1},\dots, \alpha_{l}$ linearly independent in
$\mathfrak{t}^*$, then

$$
{\rm Res}^\Lambda\left(\frac{e^{i\lambda(X)}[dX]}
{\prod_{i=1}^{l}\alpha_i(X)}\right) = \frac{1}{\det \overline{\alpha}},
$$ 
where $\overline{\alpha}$ is the 
nonsingular matrix whose columns are the 
coordinates of $\alpha_{1},\dots, \alpha_{l}$
with respect to any orthonormal basis of $\t$ defining the same orientation.

\item When $\lambda$ is of the form $\sum_{j= 1}^k  s_j \beta_j$ 
where fewer than $l$ of the $s_j$ are nonzero, then we define
the residue 
$$
{\rm Res}^\Lambda\left(\frac{e^{i\lambda(X)}[dX]}
{\prod_{i=k_1}^{k_l}\beta_i(X)}\right) 
$$ 
as the limit of the residues  at $\lambda + 
s \rho$ where $s$ is a small positive real parameter and 
$\rho \in \liet^*$ is chosen so that $\lambda + s \rho$
does not lie in any cone of dimension $l - 1$ or less spanned 
by a subset of the $\{ \beta_j\}$.
\end{enumerate}

\begin{rem}
Suppose $X^J = 1$, so all $j_m = 0 $. Then the residue 
is equal to  
$$\frac{i^{k - l}H_{\mathbf{\beta}} (\lambda)}{(2 \pi)^l}  $$
where $H_{\mathbf{\beta}}$ is the pushforward of Lebesgue measure
from $\C^k$ to $\liet^*$ 
(in other words, the Duistermaat-Heckman function for 
the linear action of $T$ on $\C^k$ via the collection of 
weights $\beta_1, \dots, \beta_k$).  (See \cite{JK:LocalizationNonabelian}
 Proposition  8.11 (ii) and \cite{GLS:Kostant}.)
 \end{rem}
 
 \begin{rem}
 By property (4), the residue equals $0$ when $\lambda = 0 $ 
provided that $k > l $, because
the functions $H_{\mathbf{\beta}}$ are piecewise polynomial
functions of degree $k - l$ and so must approach $0$ 
at $s \rho$  as $s \to 0^+$, 
when $\rho$ is a point in the cone $\Lambda^*$.
Likewise, when $k = l$ the value of the residue at $\lambda = 0$
equals the value given in (\ref{it:nonzeroresidue}). 
\end{rem}

\begin{rem}
It turns out that it is legitimate to expand the exponential in the
left hand side of (2) using a Taylor series and then to compute
the right hand side in property (2) 
 using linearity and the properties listed in (1) and (3).
 In particular, $(i\lambda(X))^m$ is a polynomial in X and can
be distributed. Then for a multi-index $I=(i_1,\dots i_l)$,
$$
\lim_{s\to 0^+}\res^\Lambda \big( \frac{X^J X^I
e^{is\lambda(X)}[dX]}{\prod_{i=1}^v \alpha_i(X)} \big) =0
$$
unless $v-(j_1+\dots j_l)-(i_1+\dots i_l)=l$.
\end{rem}

Suppose $M^T$ consists of isolated fixed points. Consider the case that
$\beta$ is of homogeneous degree, and that
$\alpha\in H^*_T(M)$ is homogeneous and complementary in degree in the sense that
$$\deg\alpha = \dim M_p-\deg\beta.$$ Then the residue theorem
states
\begin{equation}\label{eq:JKresidueisolated}
\int_{M_p}\kappa_p(\alpha\beta)= c\cdot \res^\Lambda \left(\sum_{F\in
M^T} e^{i(\Phi(F)-p)(X)}\frac{\iota^*_F (\alpha\beta)}{e_F}[dX]\right).
\end{equation}

Alternatively (see \cite{JK:LocalizationQuantization} Proposition 3.4)
the residue can be  reformulated as follows. 
For $f$ a meromorphic function of one complex variable
$z$  which is of the 
form $f(z) = g(z) e^{i \lambda z} $ where $g$ is a rational function,
we define
$$\res^+_z f(z)dz =  
\sum_{b\in\C} \res(g(z)e^{i\lambda z};z=b).$$
We extend this definition  by 
linearity to linear combinations of functions
of this form.

Viewing $f$ as a meromorphic function on the Riemann sphere and observing
that the sum of all the residues  of a meromorphic 1-form
on the Riemann sphere is $0$, we observe that 
$$ \res^+_z \left ( f(z)dz\right )  = 
- {\rm Res}_{z = \infty} \left ( f(z)dz\right ). $$

If $X  \in \liet$,
define
$$h(X)=\frac{q(X)e^{i\lambda(X)}}{\prod_{j=1}^k \beta_j(X)}$$
for some polynomial function $q(X)$ of $X$
 and some $\lambda,\beta_1,\ldots,\beta_k \in \liets$.
Suppose that $\lambda$ is
 not in any proper subspace of $\liets$ spanned by a subset of
$\{\beta_1,\ldots,\beta_k\}$. Let $\Lambda$
 be any nonempty open cone in $\liet$ contained in
some connected component of 
$$\{X\in\liet : \beta_j(X) \neq 0, 1\leq j\leq k\}.$$
Then for a generic choice of coordinate system $X=(X_1,\ldots,X_l)$ 
on $\liet$ for
which  $(0,\ldots,0,1)\in \Lambda$ we have
\begin{equation} \label{e:resdef}
\res^{\Lambda}(h(X)[dX])=  \Jac 
\res^+_{X_1} \circ  \ldots \circ \res^+_{X_l}
\left (  h(X)dX_1\ldots dX_l\right )\end{equation}
where the variables $X_1,\ldots,X_{m-1}$ are held constant while calculating
$\res^+_{X_m}$, and $\Jac$ is the determinant of any $l\times l$ matrix whose
columns are the coordinates of an orthonormal basis of $\liet$ defining the
same orientation as the chosen coordinate system. 
We assume that if $(X_1, \dots, X_l)$ is a coordinate
system for $X \in \liet$, then
$(0, 0, \dots, 1) \in \Lambda$. We  also require an additional
technical  condition on the coordinate systems, which
is valid for almost any choice of coordinate system (see
Remark 3.5 (1) from \cite{JK:LocalizationQuantization}).

One proof of (\ref{e:resdef})
 involves checking that the object on 
the right hand side of (\ref{e:resdef}) satisfies a subset of
the properties (1)-(3) above
which characterize the residue uniquely 
(see \cite{JK:LocalizationQuantization}, Proposition 3.4).

\subsection{Equivariant Morse theory}\label{se:equivariantMorsetheory}
The goal of this section is to use the restriction properties of $\alpha\in H_T^*(M)$ to find level sets on which $\alpha$ restricts to 0 under the corresponding Kirwan map. In fact we show something stronger, namely that in the case of isolated fixed points there is a basis of $H_T^*(M)$ (as an $H_T^*(pt)$-module) given by certain classes $\{\alpha_F\}_{F\in M^T_{cc}}$ indexed by the fixed point set (although the choice of these classes may not be unique). The classes $\alpha_F$ are in $I_p$ for a very specific (and generally rather large) set of regular values $p$, as we describe below. 

First consider the case that $p$ is not in the convex hull of the
image under $\Phi$ of $supp\ \alpha$. One can prove that $\alpha\in
I_p$ using equivariant Morse theory. In particular, if $p$ is not in
this convex hull, then there is some $\xi\in\t$ whose corresponding
moment map component $\Phi^\xi:= \langle \Phi,\xi\rangle$ is such that
$\Phi^\xi(\Phi^{-1}(p))<\Phi^\xi(F)$ for every $F\in supp\
\alpha$. One may choose $\xi$ so that $\Phi^\xi$ has distinct values
for distinct components of the fixed point set. Denote these components by $C_0, C_1,\dots, C_k$ and order 
them so that $i<j$ if and only if $\Phi^\xi(C_i)<\Phi^\xi(C_j)$. We
note that $\Phi^\xi$ is a Morse-Bott function. Following
\cite{AB:Yang-Mills} one defines 
\begin{align}\label{def:Mi+}
M_i^+ &:= (\Phi^\xi)^{-1}(-\infty, \Phi^\xi(C_i)+\epsilon_i)\mbox{ and}\\
\label{def:Mi-}
M_i^- &:= (\Phi^\xi)^{-1}(-\infty, \Phi^\xi(C_i)-\epsilon_i)
\end{align}
where $\epsilon_i>0$ is small enough that $C_i$ is the only critical
set in the interval
$(\Phi^\xi)^{-1}(\Phi^\xi(C_i)-\epsilon_i,\Phi^\xi(C_i)+\epsilon_i)$.
For each $i$, consider the long exact sequence in equivariant
cohomology 
\begin{equation}\label{seq:relexact}
\cdots\rightarrow H_T^*(M^+_i,M^-_i)\rightarrow H_T^*(M^+_i)\rightarrow
H_T^*(M^-_i)\rightarrow H_T^{*+1}(M^+_i,M^-_i)\rightarrow\cdots.
\end{equation}
Using the Thom isomorphism, we may identify $$H_T^*(M^+_i,M^-_i)\cong
H_T^{*-\lambda_i}(C_i),$$ where $\lambda_i$ is the index of the
negative normal bundle to $C_i$ (under $\Phi^\xi$). Atiyah and Bott
observed that the map $H_T^{*-\lambda_i}(C_i)\to H_T^*(M^+_i)$ is an
{\em injection} and hence the sequence splits 
\begin{equation}\label{seq:Cishortexact}
0\rightarrow H_G^{*-\lambda_i}(C_i)\rightarrow H_G^*(M^+_i)\rightarrow
H_G^*(M^-_i)\rightarrow 0.
\end{equation}
It follows by induction on $i$ that if $\alpha$ is 0 on $C_1$,\dots,
$C_i$, then $\alpha$ is 0 on $M^+_i$ and therefore $\alpha$  restricts
to 0 on $\Phi^{-1}(p)\subset (\Phi^\xi)^{-1}(p)$ if
$\Phi^\xi(p)<\Phi^\xi(C_{i+1})$. 

When $M$ is not K\"ahler, we may not have the familiar convexity
properties.  The convex hull of $\Phi(supp\ \alpha)$ may contain regular values that are not in the
image of the gradient flow-outs of $supp\ \alpha$: see the
non-K\"ahler example in Section~\ref{se:nonKahler}.  For these regular values, 
it is not
immediately clear that the class $\alpha$ restricts to zero at these
regular values. We use some additional Morse theory to prove
that $\alpha$ does indeed restrict to 0 on these level sets.  We first
define the {\em extended stable set} of a critical set $C$.  

\begin{definition}
Let $f$ be a Morse-Bott function on $M$, and order the critical sets
$C_0, C_1,\dots, C_k$ so that $i<j$ if and only if
$f(C_i)<f(C_j)$. Let $grad\ f$ be the gradient flow of $f$ with
respect to a compatible Riemannian metric. The {\em extended stable
set}  of $C$ is the set of points $x$ in $M$ such that there is a
sequence of critical sets $C_{i_1}, C_{i_2},\dots, C_{i_m}=C$ in which
$x$ converges to $C_{i_1}$ under the flow of $-grad\ f$, and there exist points in
the negative normal bundle of $C_{i_j}$ that converge to $C_{i_{j+1}}$
under the flow of $-grad\ f$. 
\end{definition}

The extended stable set is also called the extended flow-up (where
``up" means that the value of $f$ increases under the flow of  $grad\
f$ rather than $-grad\ f$).

\begin{rem}
Note that the extended stable set is {\em not} the same as the closure
of the stable flow-outs of the critical set $C$. A manifold with a
Morse-Bott function for which these two are not the same is the
Hirzebruch  surface $F_2$ (see Figure 1).
\end{rem}
\begin{figure}[h]
\centerline{
\epsfig{figure=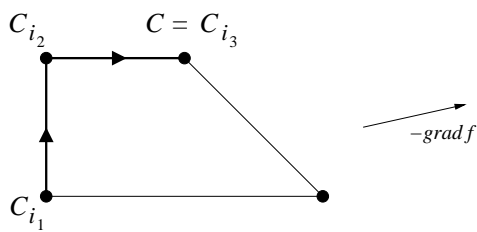,width=2.5in}
}
\smallskip
\centerline{
\parbox{4.5in}{\caption[The Hirzebr\"uch surface $F_2$]{\small The
bold lines are the image of the {\em extended stable set} of $C$ under
$\Phi$. The closure of the {\em stable flow-out} of $C$ is just the (top)
horizontal bold line.}} 
}\label{fi:Hirzebruch}
\end{figure}

A slight generalization of these methods (see
\cite{Go:Effective}\footnote{in \cite{Go:Effective} the extended
stable set was termed ``extended stable manifold,'' a misnomer since
these sets are not in general manifolds}) yields 
\begin{theorem}\label{th:Thomexistence}
For every fixed point set component $F$ there exists a homogeneous class $\alpha_F$
which has support on the extended stable set of $F$ and such that
$\iota^*_F(\alpha_F) = e_T(\nu^- F)$, the equivariant Euler class of
the negative normal bundle to $F$. We call such a class a {\em
Morse-Thom class} associated to $F$.  
\end{theorem}

It now follows easily that, in the case of isolated fixed points, any class $\alpha\in H_T^*(M)$ can be written as
 an element of the ideal generated by the classes $\alpha_F$ as $F$ ranges over the fixed point set. Let $i$
be the smallest number such that $\iota^*_{C_i}(\alpha)\neq 0$. Then
$\iota^*_{C_i}(\alpha)$ is a multiple $m$ of the equivariant Euler
class of the negative normal bundle to $C_i$. Let $\alpha_i$ be a
Morse-Thom class associated to $C_i$. Then
$\alpha-m\alpha_i$ restricts to 0 on all classes $C_j$ where $j\leq
i$. Continue inductively to express $\alpha$ in terms of the
$\alpha_i, i=1,\dots, k.$

\begin{rem}
In the case that some component $f$ of the moment map is {\em
Palais-Smale}, i.e. the (ordinary, non-extended) unstable and stable
manifolds intersect transversally, there is a unique Morse-Thom class for each component of the fixed point set:
Suppose $F$ and $G$ are two isolated fixed points and that the
unstable manifold of $F$ intersects the stable manifold of $G$
transversally and non-trivially. If $F\neq G$, then the intersection
is not zero-dimensional and $\dim \nu^+_G+\dim \nu^-_F-\dim M > 0$,
where $\nu^+_G$ is the positive normal bundle of $G$ in $M$ and
$\nu^-_F$ is the negative normal bundle of $F$ in $M$, both with
respect to $f$. Since $\dim \nu^+_G=\dim M-\dim\nu^-_G$, it follows
that $\dim\nu^-_F>\dim\nu^-_G$.  The degree of any Morse-Thom
class out of $G$ is equal to  
 $\dim \nu^-_G$, and thus homogeneity assures that the class
associated to $G$ is unique.  This is the case for coadjoint orbits of
complex reductive groups. A similar argument works in the case that
$F$ and $G$ are not isolated. 
\end{rem}

Let $p$ be a regular value of $\Phi$, $f$ a Morse-Bott function on $M$, and $C_0,\dots, C_k$ the connected components of the critical set. We choose $f$ generically enough so that we may order the critical sets  by $f(C_i)<f(C_j)$ for $i<j$. Suppose that 
$\alpha\in H_T^*(M)$ is an equivariant Euler class associated to the connected component $C_i$ of the critical set. 
Suppose that
$\Phi^{-1}(p)$ does not intersect the extended stable set associated to $F$. Then $\alpha|_{\Phi^{-1}(p)}=0$, for
$\Phi^{-1}(p)$ is a subset of a space equivariantly homotopic to $M^-_i$ (see \cite{Go:Effective}), and $\alpha|_{M_i^-}=0$ by assumption that it is an Euler class associated to $C_i$. 
Thus $\alpha$ restricted to a level set $\Phi^{-1}(p)$
is non-zero only if this level set is in the extended stable set of a
fixed point component $F\in supp\ \alpha$. 
We have proven
\begin{theorem}
Let $\alpha\in H^*_T(M)$, where $M$ is a compact Hamiltonian $T$-space
with moment map $\Phi$. Suppose $p$ is a regular value of $\Phi$, and
that for some component $\Phi^\xi$ of $\Phi$,  $\Phi^{-1}(p)$ does not
intersect the extended stable set with respect to $\Phi^\xi$ of any of the points in the support of
$\alpha$. Then $\alpha|_{\Phi^{-1}(p)}=0$.
\end{theorem}
\noindent We will need this theorem to prove Theorem~\ref{th:main} in the general case, in which we are not guaranteed that the supports of Morse-Thom classes are convex. We note, however, that we have proven
\begin{corollary}\label{cor:halfspacesupport}
Let $H\subset \t^*$ be a rational hyperplane passing through $p$. We
orient $H$ by realizing it as a level set of a vector $\xi\in\t$ and
denoting by $H^+$ the set of all points $q\in\t^*$ such that
$\xi(q)>\xi(H)$, and by $H^-$ the set of all points $q\in\t^*$ such
that $\xi(q)<\xi(H)$. If $\Phi(supp\ \alpha)\subset H^+$ or
$\Phi(supp\ \alpha)\subset H^-$, then $\alpha\in I_p$. 
\end{corollary}

 A much harder fact to prove is that the set of these classes, each of
which has this vanishing property for some hyperplane $H$, generate
the kernel of the Kirwan map. This is the content of the
Tolman-Weitsman theorem \cite{TW:abelianreduction}.

Notice that regular values are never on walls of the moment
polytope. If $p$ is regular, then $\Phi^{-1}(p)$ has a locally free
$T$ action, and in particular points in the level set can be fixed  by
at most 0-dimensional subsets of $T$.

\section{The Main Theorem and its Proof}

We begin with a proposition about the behavior of $I_p$ as $p$ moves
within a chamber. 

\begin{prop}\label{pr:samechamber}
Let $M$ be a compact Hamiltonian $T$ space with moment map $\Phi$. If
$p$ and $q$ are two regular values of $\Phi$ in the same chamber of
the moment polytope, then $I_p=I_q$. 
\end{prop}

\begin{proof}
Since any two points $p$ and $q$ in the chamber $C$
may be connected by a piecewise linear
path lying entirely in $C$, the Proposition
follows from the following theorem
 of Guillemin-Kalkman \cite{GK:residue} and Martin
\cite{Ma:wallcrossing}.

\begin{theorem}[\cite{GK:residue},\cite{Ma:wallcrossing}] Let there be
a piecewise linear path of regular values between two regular values
$p$ and $q$ of the moment map. Let $\eta \in H^*_T(M)$. 
Then for any $\zeta \in H^*_T(M)$ we have
$$\int_{M_{p}} \kappa_{p} (\zeta) \kappa_{p}(\eta)
=   \int_{M_{q}} \kappa_{q} (\zeta) \kappa_{q}(\eta). $$
\end{theorem}

Since $\kappa_{p}$ and $\kappa_{q}$
are surjective, it follows that
$\kappa_{p}(\eta) = 0 $ if and only if
$\kappa_{q}(\eta) = 0$
(using Poincar\'e duality). \end{proof}

We now assume that $T$ acts on $M$ with isolated fixed points, and prove the harder part of Theorem \ref{th:main}, namely
\begin{prop}
Let $M$ be a compact Hamiltonian $T$ space with isolated fixed points
and moment map $\Phi$. If $p$ and $q$ are two regular values of $\Phi$
in distinct chambers 
of the moment polytope, then $I_p\not=I_q$. 
\end{prop}

We prove this by construction. For any $p$ and $q$ in distinct
chambers $C_p$ and $C_q$ respectively, we find a  class $\gamma$ with certain restriction
properties that make computations easy. We then show that the
Jeffrey-Kirwan localization theorem implies that $\gamma\in I_p$ but
$\gamma\not\in I_q$. 

We note that any wall can be oriented as follows. Suppose $W\subset
\t^*$ is a codimension-1 wall of $\Delta$. Choose $\xi\in\t$
perpendicular to $W$ in the sense that $W$ lies in the intersection of
$\Delta$ with a level set of $\xi$.  We write $W^+$ for $\{x\in \t^* \
|\ \xi(x)>\xi(W)\}$ and $W^-$ for $\{x\in \t^*\ |\
\xi(x)<\xi(W)\}$. Thus $W^+$ and $W^-$ are open half spaces in
$\t^*$. We denote by $H$ the hyperplane of $\t^*$ containing $W$. Define $H^+=W^+$ and $H^-=W^-$. 

Let $p$ and $q$ be regular values of the moment map in distinct
connected components of $\Delta$. Choose $W$ a codimension-1 wall of
the chamber of $p$,  and orient $W$ so that $p\in W^+$ and
$q\in W^-$. Note that $W$ need not equal $C_p\cap W$. For every fixed point $F$, let $\{\beta^F_i\}$ denote the
set of weights of the $T$ action on $T_FM$.

We begin by finding a special subset $P$ of the fixed point set
$M^T$. For every $F$, let $J_F$ be the index set with $i\in J_F$ if
and only if $\Phi(F)+\beta_i^F\in H^+\cup H$. 
The set $P$ consists of the fixed points $F$ that satisfy two conditions:
\begin{enumerate}
\item $\Phi(F)\in W$, and
\item $p$ is in the positive linear span of
$\{\Phi(F)+\beta^F_i\}_{i\in J_F}$. We will write $p\in
\Phi(F)+\langle\beta_i^F\rangle_{i\in J_F}^+$. 
\end{enumerate}

\begin{lemma}
$P$ is non-empty.
\end{lemma}
\begin{proof}

Let $W$ be a wall adjacent to the chamber $C_p$ containing $p$, and $H$
the
hyperplane containing $W$. The wall $W\subseteq H$ is contained in
the
image under $\Phi$ of a connected component $M_0$ of the fixed point set $M^S$
of a
circle $S\subseteq T$. This component $M_0$ is itself a
compact symplectic manifold with Hamiltonian $T/S$ action, and must
therefore have fixed points.  Therefore, there are fixed points that
satisfy condition (1).

Since $W$ is a wall on the chamber of $p$, we may assume that $p$ is
very
close to $W$. As above, $M_0$ is a connected component of $M^S$
such that
$W\subseteq \Phi(M_0)\subset H$.  Consider the orthogonal
projection $\pi:\mathfrak{t}^*\to H$.  Since $p$ is close to $W$, we
may
assume that $\pi(p)\in W$.

We note that the maximal value of $\xi \circ \Phi$ is not attained on $M_0$ since $W$ is necessarily an internal wall.
Therefore, at every fixed point $F\in M_0^T$, the action of $T$ on
$T_FM_0$ must have at least one weight pointing into $H^+=W^+$.

Now consider the weights pointing along $H$. The $T/S$ action on $M_0$
is effective, and thus the image of the moment map for this action is
a convex set in $Lie(T/S)^*$, identified with a subset of $H$
containing $W$. It follows that, for any point in $W$ (and in
particular, for $\pi(p)$), there is a fixed point  $F\in M_0^T$ whose
weights have positive span containing this point.  
Choosing a fixed point $F$ with $\pi(p)$ contained in the positive
span of the weights $\{\beta_i^F\}$ at $F$ pointing into $H$, we may
move $p$ close enough to $W$ so that $p$ is contained in the positive
span of $\{\Phi(F)+\beta^F_i\}$.

This proves that there is an $F\in M^T$ satisfying the two conditions
in the definition of $P$.
\end{proof}

We now choose an appropriate dual cone $\Lambda^*$ in $\t^*$. Consider
the set $\{\beta_i^F\}_{i,F}$ of all weights that occur at any fixed
point. Choose any $F_0\in P$. Since $p\in \Phi(F_0)+\langle
\beta_i^{F_0}\rangle_{i\in J_{F_0}}^+$, there exists a subset of
weights $$\{ \beta_i^{F_0}\}_{i\in {I_{F_0}}}\subset \{
\beta_i^{F_0}\}_{i\in J_{F_0}}$$ indexed by $I_{F_0}\subset J_{F_0}$
such that 
\begin{enumerate}
\item[(a)] $\langle \beta_i^{F_0}\rangle^+_{i\in I_{F_0}}$ does not
contain a line, and  
\item[(b)] $p\in \Phi(F_0)+\langle \beta_i^{F_0}\rangle_{i\in I_{F_0}}^+$.
\end{enumerate}
Choose any $X\in\t$ with $\beta_i^F(X)\neq 0$ for all $i, F$ and satisfying:
\begin{enumerate}
\item $\beta_i^{F_0}(X)<0$ for all $i\in I_{F_0}$, and
\item\label{it:noH-vectors} If $\Phi(F)\in H$ and $\Phi(F)+\beta_i^F
\in H^-$, then $\beta_i^F(X)>0$. 
\end{enumerate}
The vector $X$ defines a polarization of the weights.
Let $\gamma_i^F = \beta_i^F$ if $\beta_i^F(X)>0$ and $\gamma_i^F =
-\beta_i^F$ if $\beta_i^F(X)<0$. Finally, let $\Lambda^*$ be the cone
generated by the span of the vectors $\{\gamma_i^F\}_{i,F}.$

Armed with $\Lambda^*$, we are now prepared to find a class $\gamma\in
I_q$ with $\gamma\not\in I_p$. We first show that the choice of
$\Lambda^*$ ensures that certain fixed points (those whose image under
$\Phi$ lie in $H^+$) will not contribute to the residue for any
$\beta\in H_T^*(M)$. By Proposition \ref{pr:samechamber}, we may
assume that $p$ is close to $H$. The choice of polarization implies that (for $p$ close enough to $H$) if $\Phi(F)\in H^+$, then $\Phi(F)\not \in p+\Lambda^*$.  In
particular by Property (1)(c) of the residue, 
$$
\res^\Lambda \frac{e^{i(\Phi(F)-p)}}{\prod_{i=1}^l\alpha_i}=0
$$ 
for {\em any}  $\alpha_1,\dots, \alpha_l\in \Lambda^*$, whenever
$\Phi(F)\in H^+$. 

We now show that there is a certain set of classes $\alpha\in
H_T^*(M)$ which restrict to 0 on fixed points $F$ with $\Phi(F)\in
H^-$. Recall that $H$ is level set of $\xi\in \t$. Consider $\eta\in
\t$, a small perturbation of $\xi$ with the properties that 
\begin{enumerate}
\item For $F_1, F_2\in M^T$ with $\Phi(F_1), \Phi(F_2) \in H$ but
$\Phi(F_1)\neq \Phi(F_2)$, we have $\eta(\Phi(F_1))\neq
\eta(\Phi(F_2))$. 
\item For any $F\in M^T$ with $\Phi(F)\in H^+$ and any $F_1\in M^T$
with $\Phi(F_1)\in H$, we have $\eta(\Phi(F_1)
)<\eta(p)<\eta(\Phi(F)).$ 
\end{enumerate}
 
 Consider $\eta\circ\Phi$ restricted to $P$. Let $F_1\in P$ be the
fixed point where $\eta\circ\Phi|_P$ is maximized. Recall that by
construction, $\Phi(F_1)\in H$.  
 By Theorem \ref{th:Thomexistence} there exists a Morse-Thom class $\alpha=\alpha_{F_1}$ which has support (in the sense of definition \ref{def:support}) on the extended flow-up out of $F_1$ along $grad\ \eta$. In particular, 
$$\alpha|_{F}=0 \hspace{.2in}\mbox{if}\hspace{.2in} \Phi(F)\in H^-.$$
It now follows immediately that
\begin{lemma}\label{le:relevantEulerclass} Let $\alpha=\alpha_{F_1}$ as above. Then
$\beta\alpha\in I_q$ for any $\beta\in H_T^*(M)$.
\end{lemma}
\begin{proof}
Using the cone $\Lambda^*$, the only terms that contribute to the
integral over $M_q$ are fixed points $F$ such that $\Phi(F)-q\in \Lambda^*$, or $q\in \Phi(F)-\Lambda^*$.
Since $q\in H^-$, it follows that the only possibly contributing $F$ are those whose images under the moment map
are in $H^-$. However, $\beta\alpha$ restricts to 0 on all these
points (as $\alpha$ does), which implies the residue contributions are
0. 
\end{proof}

It is left to show there exists $\beta\in H_T^*(M)$ such that
$\int_{M_p}\kappa_p(\beta\alpha)\neq 0$.  We first simplify the
integral by showing that there is only one possibly nonzero term in
the residue formula (\ref{eq:JKresidueisolated}). 
\begin{claim} Let $\alpha$ be a Morse-Thom class of $F_1$ and
$\beta\in H_T^*(M)$ be any class.  
Then
\begin{equation}\label{eq:JKsingleterm}
\int_{M_p}\kappa_p(\beta\alpha)= c\cdot \res^\Lambda
e^{i(\Phi(F_1)-p)(X)}\frac{\iota^*_{F_1}(\beta\alpha)}{e_{F_1}}[dX]. 
\end{equation}
where $e_{F_1}$ is the equivariant Euler class at $F_1$.
\end{claim}

\begin{proof}
Let $F$ be a fixed point with $\Phi(F) \in H^+$. Then we showed that
$p\not\in\Phi(F)-\Lambda^*$, and the corresponding residue term in
Equation (\ref{eq:JKresidueisolated}) is therefore 0. 

Suppose instead that $F$ is a fixed point with $\Phi(F)\in H^-$. Then
$\alpha|_F=0$ implies that any corresponding residue term is 0. 

Finally, suppose that $F$ is a fixed point with $\Phi(F)\in H$. If
$F\not\in P$, then by construction $p\not\in \Phi(F)-\Lambda^*$ and
therefore the corresponding residue is 0. On the other hand, if $F\in
P$ and $F\neq F_1$, then $\alpha|_F=0$ since $\alpha$ is supported on
the flow-up out of $F_1$ along $grad\ \eta$ and $F_1$ is where
$\eta\circ\Phi$ attains its maximum among points in $P$. Therefore the
only remaining possible contribution to the integral
(\ref{eq:JKsingleterm}) is from $F_1$. \end{proof} 

The final step is to prove our original claim, namely that 
\begin{lemma} There exists $\beta\in H_T^*(M)$ such that
$\beta\alpha \not\in I_p$, where $\alpha$ is a Morse-Thom
class associated to $F_1$, as above. 
\end{lemma}
\begin{proof}
We show that, for an appropriate choice of $\beta$, the single term in
\eqref{eq:JKsingleterm} is nonzero.  

Let $R$ be the set of all weights of the $T$ action at $F_1$. Let
$R^+=\Lambda^*\cap R$ and $R^-=R\backslash R^+$.  Since $F_1$ is
isolated, the normal bundle is topologically 
trivial (although not equivariantly) and so the Euler class of the
normal bundle and that of the negative normal bundle with respect to
$\eta\circ\Phi$ are 
products of weights. In particular
$$
e_{F_1}=\prod_{w\in R} w
$$
and
$$
\iota^*_{F_1}(\alpha) = \alpha|_{F_1}= \prod_{w\in R^-} w.
$$
We simplify
\begin{equation}\label{eq:singletermroots}
\frac{\iota^*_{F_1}(\beta\alpha)}{e_{F_1}} = \frac{\iota^*_{F_1} (\beta)
\prod_{w\in R-}w}{\prod_{w\in R}w} = \frac{\iota^*_{F_1}
(\beta)}{\prod_{w\in R^+} w} 
\end{equation}
We choose $\beta$
carefully. Let
$\gamma_1,\dots,\gamma_l$ be any basis of $\t^*$ chosen among elements
of $R^+$. Let  
$$
\beta=\prod_{w\in R^+-\{\gamma_1,\dots,\gamma_l\}}w.
$$ 
As $\beta$ is in the image of $H^*_T(pt)$, $\iota^*_F\beta=\beta$ for
all fixed points $F$. Thus the formula~(\ref{eq:singletermroots})
becomes 
$$
\frac{\iota^*_{F_1} (\beta)}{\prod_{w\in R^+}
w}=\frac{1}{\prod_{i=1}^l\gamma_i}. 
$$
By Property \ref{it:nonzeroresidue} of the residue formula, the residue 
$$
c\cdot \res^\Lambda
\left(e^{i(\Phi(F)-p)(X)}\frac{1}{\prod_l\gamma_i}[dX]\right)= 
c\cdot\frac{1}{\det \overline{\gamma}}
$$ where
$\overline{\gamma}$ is a matrix whose columns are $\gamma_i$ in an
appropriate basis. Since the $\gamma_i$ are linearly independent and $c$
is non-zero, the residue is non-zero. This is the only contribution to
the sum~(\ref{eq:JKresidueisolated}), so we have shown the integral of
$\beta\alpha$ over the reduced space at $p$ is non-zero. In other
words, $\gamma=\beta\alpha\not\in I_p$.  In particular, this implies
that $\alpha\not\in I_p$, as well as $\beta\not\in I_p$.  This
completes the proof of the lemma. 
\end{proof}

\section{Applications of the main theorem}

\subsection{Schubert classes and counting chambers}\label{se:schubert}
In the case that $M$ is a coadjoint orbit of any complex reductive
group, a specific basis for $I_p$ has been described in terms of the
choice of symplectic structure on the coadjoint orbit and the
reduction point $p$ (\cite{Go:Weightvarieties},
\cite{GM:coadjointreductions}). The generators are {\em permuted
Schubert classes} defined by their duality properties with certain
subvarieties of $M$.  They can also be uniquely described by their
restriction properties due to Kirwan injectivity. We present the
latter description here. Theorem~\ref{th:Schubertrestriction} can be
taken as a definition for readers unfamiliar with Schubert classes. 

 Let $G$ be a compact semi-simple Lie group with Borel subgroup $B$
and maximal torus $T$. Let $M$ be a coadjoint orbit of $G$ through
$\lambda\in\t^*$. We identify the Weyl group $W:=N(T)/T$ with the
fixed point set by associating $\lambda$ to the trivial coset, and the
remaining points can be identified by the transitive action of $W$ on
the fixed points. We write $\lambda_\sigma$ to indicate the fixed
point corresponding to $\sigma$. Let $f$ be a moment map for a generic
$S^1$ action on $M$, such that $f|_{M^T}$ is minimized at $\lambda$. 
\begin{theorem}\label{th:Schubertrestriction} For every $\sigma\in W$
the associated Schubert class $\beta_\sigma\in H_T^*(M)$ is the unique
homogeneous class defined by the following restriction properties. 
\begin{enumerate}
\item $\beta_\sigma|_{\lambda_\sigma} = e(\nu^-_f \lambda_\sigma)$,
the equivariant Euler class of the negative normal bundle (with
respect to $f$) of the fixed point associated to $\sigma$. 
\item $\beta_\sigma|_{\lambda_\tau} =0$ unless $\tau\geq\sigma$ in the
Bruhat order of $W$. 
\end{enumerate}
\end{theorem}
Schubert classes form a basis for $H^*_T(M)$ as a module over
$H^*_T(pt)$. If we use the same construction while allowing $f$ to be
minimized at other fixed points, we obtain the {\em permuted Schubert
classes} 
$$\beta_\sigma^\tau := \tau\cdot \beta_{\tau^{-1}\sigma},$$
whose restriction properties can be described in a similar fashion.

Permuted Schubert classes are precisely the classes predicted by
Theorem \ref{th:Thomexistence}. In the case of coadjoint orbits the
{\em extended} flow-up is precisely the closure of the flow-up from
any fixed point. They form (permuted) Schubert varieties. 

There is a very simple description of ideals $I_p$ for coadjoint
orbits of complex reductive groups \cite{Go:Weightvarieties},
\cite{GM:coadjointreductions}. For any regular value $p$, the ideal
$I_p$ is generated by the permuted Schubert classes
$\beta_\sigma^\tau$ such that $p$ is not in the convex hull of $supp\
\beta_\sigma^\tau$. Since the number of chambers in $\Delta$ is
determined by the number of distinct ideals $I_p$, this simple
description offers the possibility of an algebraic method for
ascertaining the number of chambers in the moment polytope for
coadjoint orbits $SU(n)$, for example. Of course, these numbers depend
on the non-zero value $\lambda\in\mathfrak{t}^*_+$ in the positive
Weyl chamber, whose coadjoint orbit we are considering. For $SU(3)$,
the number of chambers is six or seven, and does not offer much of a
challenge. E. Rassart has computed this number for $SU(4)$ orbits
\cite{Ra:numberchambers} using geometric arguments, and has found
numbers in the hundreds (depending on the value of $\lambda$), but the
$SU(5)$ case is already out of reach. It would be interesting to find
a bound on the number of chambers for general $n$. 

\subsection{Relation to a wall-crossing formula}

There are several places in the literature in which the behavior of
the reduced space, or integrals over these spaces, is studied as the
regular value moves over exactly one codimension-1 wall of the moment
map. Most notable, perhaps, is the formula due to Guillemin-Kalkman
\cite{GK:residue}, and independently S. Martin
\cite{Ma:wallcrossing}. We present here a short exposition on this
result, and then show that Theorem \ref{th:main} allows us to gain
some information about the reduction at singular values of the moment
map in the case that the fixed points are isolated. 

Let $L$ be a hyperplane of critical values of $\Phi$
such that the line segment between two regular values
$p$ and $q$ intersects $L$ but intersects no
other hyperplanes of critical values of $\Phi$.
Let $N$ be a component of the fixed point set of a
circle $S \subset T$
such that $\Phi(N) \subset L$. Then
Guillemin-Kalkman give
$$\int_{M_{p}} \kappa_{p} (\zeta \eta)
- \int_{M_{q}} \kappa_{q} (\zeta \eta) =
\int_{N_{\rm red}} \kappa_{T/S} \left ( {\rm Res}
\left [ \frac{\iota_N^* (\zeta \eta)}{e_N} \right ] \right ). $$
Here, $\Phi_S (M_p) > \Phi_S (M_q)$, where $\Phi_S$ is the
component of $\Phi$ corresponding to  $S$.
This requires the choice of a basis element in the Lie algebra of $S$,
and the same basis element must be used to define the residue. 
Note that $N_{\rm red}$ is the reduced space of $N$
with respect to the action of $T/S$,
and $e_N$ is the equivariant Euler class of the normal
bundle of $N$ in $M$; $\kappa_{T/S} $ is
the Kirwan map for the action of   $T/S$ on $N$.

The operation ${\rm Res}$ is defined as follows.
Introduce a basis element $X$ for ${\rm Lie}(S)$
and basis elements $Y_1, \dots, Y_n$ for
${\rm Lie }(H)$, where $H$ is a codimension-$1$ subtorus transverse to $S$, so $H \cong T/S$.
The normal bundle $\nu$ to $N$ can be assumed
to be a sum of line bundles:
$$\nu = L_1 \oplus \dots \oplus L_l.$$
We put
$$ c_1^T(L_i) = m_i X  + \mu_i + \sum_j f_i^j Y_j$$
where $m_i \in {\bf Z}$ is the weight of the
representation of $S$ on $L_i$ and
$\mu_i + \sum_j f_i^j Y_j$ is a Cartan representative
of $c_1^H(L_i)$
(so $\mu_i \in \Omega^2(N)$ and $f_i^j \in \Omega^0(N)$).
Write
$$\frac{1}{c_1^T(L_i) } = \frac{1}{m_i X(1+ \sum_j f_i^j Y_j/X
+ \mu_i/X)}$$
and expand
$$ \frac{1}{1+A} = 1-A+A^2 - \dots $$
using $A = (\sum_j f_i^j Y_j + \mu_i)/X. $
By expanding into a product of power series in
$X$ we can rewrite $\iota^*_N(\zeta \eta)/e_N$
as a sum $\sum_{j = - \infty}^{m_0} \beta_j X^j $
where each $\beta_j$ is a polynomial in
$Y_1, \dots, Y_n$ with coefficients in $\Omega^*(N)$. In particular,
$\beta_j\in \Omega_{T/S}^*(N)$. 
We define
$$ {\rm Res}(\iota^*_N (\zeta \eta)/e_N) = \beta_{-1}. $$
It is shown in \cite{GK:residue} that $\beta_{-1}$ is well defined 
independent of the choice of the variables $X$ and $Y_j$.

\begin{rem} \label{rem4.2}
This is the residue at $X = \infty$
of the meromorphic function of $X$ defined
by $\iota^*_N (\zeta \eta)/e_N$. Thus it equals minus the
sum of the residues at finite values of $X$
(using the fact that the sum of the residues of a meromorphic 1-form over the
Riemann sphere is zero). Here the meromorphic 1-form depends on one
complex variable $X $ in the Riemann sphere, 
but takes
values in $\Omega_{T/S}^*(N)$ 
\end{rem}

The technique for crossing one wall can of course be inductively
applied to more than one wall, obtaining a formula for the difference
in the integral over reduced spaces at {\em any} two regular
values. However the resulting formula does not imply that the chambers
are distinguished by the ideals $I_p$, and indeed it may be a lengthy
computation using this method to obtain that this difference is
non-zero for some class.

We finish this section with a small result on the behavior of the
classes that distinguish ideals. Let $p$ and $q$ be regular values
with a line between them crossing exactly one singular value $q_1$ of
$\Phi$. By Theorem \ref{th:main}, there is a class $\alpha$ such that
$\alpha\in I_q$ but $\alpha\not\in I_p$. We may choose a class with
the additional property that $\deg\alpha=\dim M_p$. By the
wall-crossing formula, 
\begin{equation}\label{eq:wall}
\int_{M_{p}} \kappa_{p} (\alpha)-0 =
\int_{N_{\rm red}} \kappa_{T/S} \left ( {\rm Res}
\left [ \frac{\iota_N^* (\alpha)}{e_N} \right ] \right )\neq 0. 
\end{equation}
Choose a  $T$-invariant Riemannian metric on $M$.
Then Equation (\ref{eq:wall}) implies that the image of $\alpha$ under
the restriction 

$$
\xymatrix{
& {H^*_T(M)}\ar[r]^{i^*} & {H^*_T(N)} \ar[r]^{p^*} & H^*_T(S(\nu N))
\ar `r[d] `[l]`[llld]`[d]_{\cong}[dll] &  \\
 &   H^*_{T/S}(S(\nu N)/S) \ar[r]^{\hspace{.3in}\pi_*} &
H^{*-k}_{T/S}(N) \ar[r]^{\kappa_{T/S}} & H^{*-k}(N_{red}) 
}
$$
is nonzero. The above composition is the ``localization map" defined
by \cite{Ma:wallcrossing}, Definition 5.1. Here $S(\nu N)$ is the unit
sphere bundle (with respect to the invariant metric) of the normal
bundle to $N$, and $N_{red}$ is the reduced space at the singular
value $q_1$. The maps $i^*$ and $p^*$ are induced by the inclusion
$i:N\hookrightarrow M$ and the projection $p:S(\nu N)\to N$,
respectively. The isomorphism on the curved arrow follows from the
fact that $S$ acts locally freely on the sphere bundle (and the
cohomology is in rational coefficients). Lastly, $\pi_*$ is
``integration over the fibers", and reduces degree by $k$, the
dimension of the fibers of the (weighted) projective bundle $S(\nu
N)/S\to N$, according to the induced orientation defined in
\cite{Ma:wallcrossing}. In particular, not only are the fundamental
classes described by Theorem~\ref{th:Thomexistence} nonzero when
restricted to level sets near $C$ and in the flow-up from $C$, but
these classes are nonzero when restricted to the {\em singular walls}
on the boundary (near $C$) of the cone out of $C$.

It is shown in Section 4 of \cite{GK:residue} how to use the formula
(\ref{eq:wall})
inductively to give a formula 
for $\int_{M_{p}} \kappa_p (\alpha)$ in terms of the restrictions of 
$\alpha$ to the components of the fixed point set of $T$.
This formula is obtained by drawing a line $l$ from $p$ to the 
boundary of $\Phi(M)$ which avoids all codimension-2 walls.
Next,  for each intersection $p_j$ of $l$ with a codimension-1
wall $W_j$, one draws a line $l_j$ within $W_j$ from $p_j$ to 
the boundary of $\Phi(M) \cap W_j$. This process gives rise to 
a graph (the term ``dendrite'' is used in \cite{GK:residue})
originating at $p$ 
and terminating in points $\Phi(F)$ where $F$ is a fixed point of 
$T$ on $M$. The quantity $\int_{M_p} \kappa_p(\alpha)$ may
  thus be expressed in terms of the data
  $\alpha|_F$ for a distinguished collection of
  fixed points $F$.

One can easily deduce from the description of the residue given in 
(\ref{e:resdef}) that Guillemin and Kalkman's iterated residue
is equivalent to Jeffrey and Kirwan's.\footnote{Both the 
Jeffrey-Kirwan residue and the Guillemin-Kalkman residue
depend on a set of choices: for Jeffrey-Kirwan one chooses
the cone $\Lambda$, while for Guillemin-Kalkman one chooses a ``dendrite''.} See also Kalkman's proof
of their equivalence \cite{Ka:res}.

\subsection{A non-K\"ahler example}\label{se:nonKahler}

We now explore the main theorem in the context of a non-K\"ahler
manifold $M$, introduced by Sue Tolman \cite{To:nonkahler}.  In this
non-K\"ahler setting, the extended stable sets are not manifolds, or
even complex varieties. We find Morse-Thom classes on an
unusual extended flow, whose image under $\Phi$ is not convex. 

\begin{center}
\begin{figure}[h]
\epsfig{figure=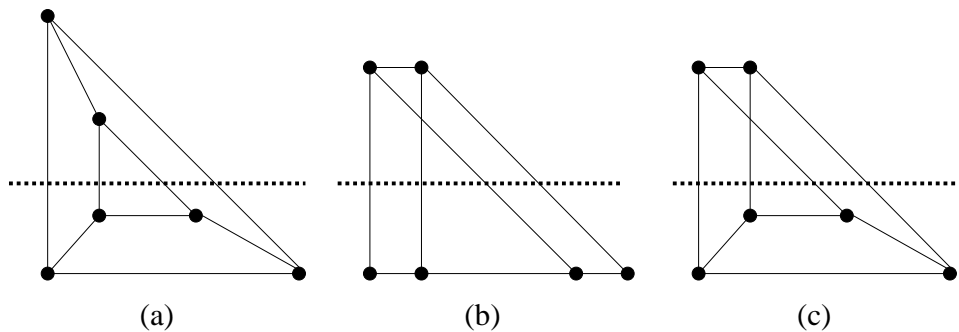,width=5in}
\caption{The non-K\"ahler manifold $M$ with
moment polytope in (c) is achieved as a symplectic gluing of the
bottom half of (a) to the top half of (b).}\label{fig:cut} 
\end{figure}
\end{center}

We first sketch Tolman's construction of $M$.  It is obtained from two
manifolds, each diffeomorphic to $N = \C P^2\times \C P^1$.  This
manifold $N$ is a toric variety with an effective Hamiltonian $T^3$
action.  We restrict our attention to two different $T^2$ actions on
$N$, whose moment polytopes are given in Figure~\ref{fig:cut}~(a) and
(b).  We then take a symplectic cut of each of these manifolds, along
the dotted line shown in the figure.  Notice that the symplectic
slices, whose moment images are the intersection of the dotted line
with the moment polytopes, are equivariantly symplectomorphic. We use the technique of symplectic gluing to glue the bottom half
of Figure~\ref{fig:cut}~(a) and the top half of
Figure~\ref{fig:cut}~(b) to construct $M$, whose moment polytope is
shown in Figure~\ref{fig:cut}~(c). 

There is a Hamiltonian $T^2$ action on $M$, with isolated fixed
points.  The moment polytope is given in Figure~\ref{fig:cut}~(c), and
the vertices of this graph correspond to the fixed points. The
restriction of a cohomology class $\alpha\in H_T^*(M)$ to a fixed
point $F\in M^T$ is an element of $H_T^*(pt)\iso \C[x,y]$ (see
Theorem~\ref{th:injectivity}).  Thus, we may describe a class in the
equivariant cohomology of $M$ by giving a polynomial in two variables
at each vertex of the moment polytope.

\begin{center}
\begin{figure}[h]
\epsfig{figure=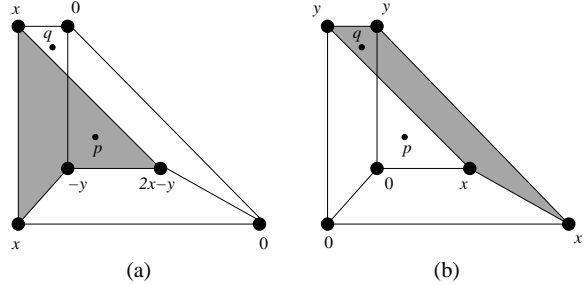,width=3.0in}
\caption{Figure (a) shows a class $\alpha \in I_q$ such that $\alpha\not\in I_p$, and (b) shows
a class in $\beta\in I_p$ such that  $\beta\not\in I_q$.}\label{fig:classes} 
\end{figure}
\end{center}

In Figure~\ref{fig:classes}, we label two points $p$ and $q$ which are
regular values of the moment map for $T$ acting on $M$.  They are in
different chambers of the moment polytope, and hence by the main
theorem, the kernel of the Kirwan map is different for each of these
values. We may apply the techniques from the proof of the main theorem
to construct equivariant cohomology classes distinguishing these
chambers.   In  Figure~\ref{fig:classes}~(a), we demonstrate a class $\alpha$
in $I_q$ but not in $I_p$, and in (b), we show a class $\beta$ in $I_p$, not
in $I_q$.  These classes distinguish these two ideals. 

The shaded region represents the image under the moment map of the
flow-out of the support of each class, for a given choice of $f$.  The
fact that the region in (a) is non-convex reflects the fact that 
the manifold $M$ is non-K\"ahlerizable. We may contrast
this to the setting given in Section~\ref{se:schubert}, where the
classes constructed that distinguish among ideals are supported on
permuted Schubert varieties. While these varieties may not be
manifolds, they are complex, and hence have convex moment image.

\begin{center}
\begin{figure}[h]
\epsfig{figure=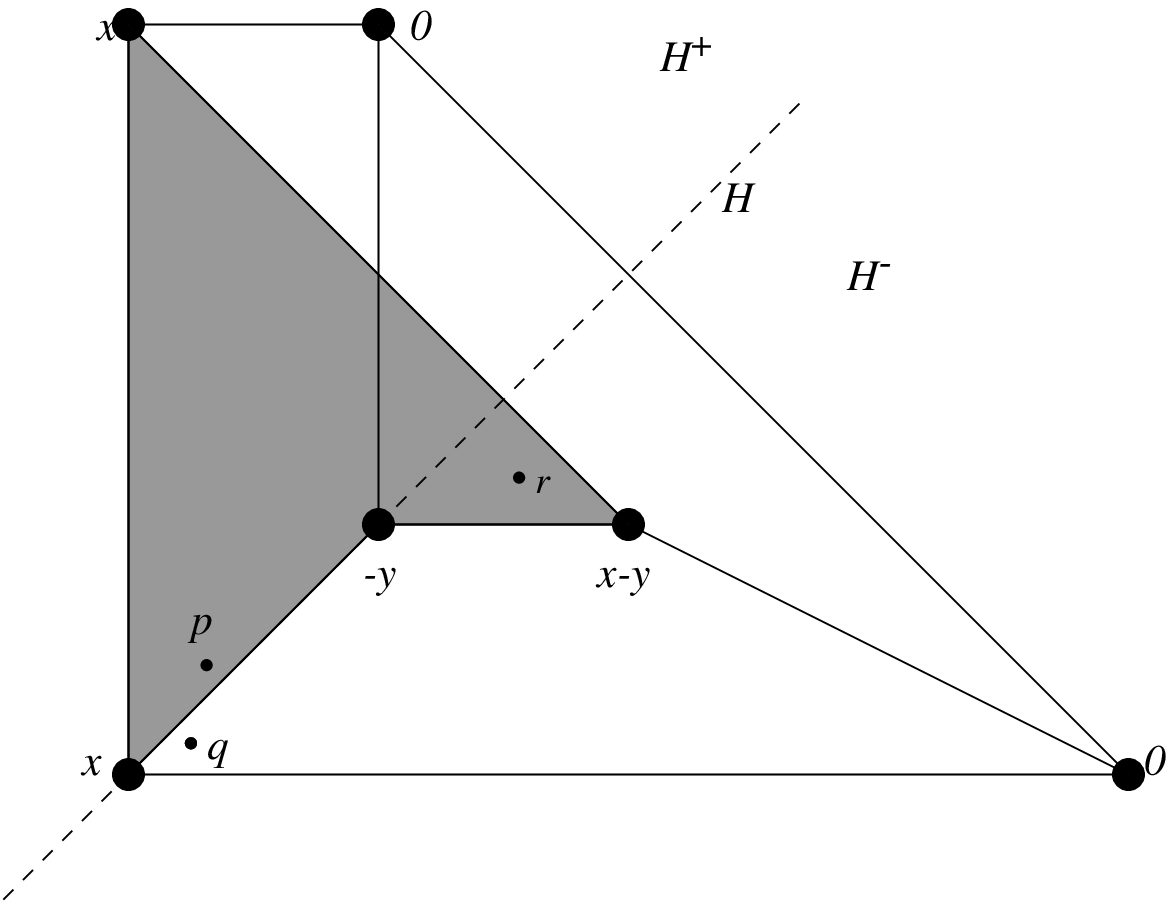,width=2.3in}
\caption{}\label{fig:nonconvex} 
\end{figure}
\end{center}

One important point about these extended flows is the following. Suppose $p$ and $q$ are in distinct chambers, and a wall $W$ adjacent to the chamber $C_p$ separates them, such that $p\in H^+$ and $q\in H^-$, where $H$ is the hyperplane containing $W$. Let the class $\alpha$ be a Morse-Thom class (associated to a point on the wall $W$) constructed by Lemma~\ref{le:relevantEulerclass}, so that $\alpha\in I_q$ but not in $I_p$. When the extended stable set is convex, then Corollary~\ref{cor:halfspacesupport} implies that $\alpha\in I_r$ for any $r\in H^-$ (the negative half space). However, without this convexity requirement, $\alpha$
 is not necessarily in $I_r$ for all $r \in H^-$. Let $p, q$ and $r$ be as indicated in Figure~\ref{fig:nonconvex}. Let $W$ be the wall slightly to the left of $q$, as indicated. The class $\alpha$ supported on the shaded region will be in $I_q$ and not in $I_p$, by the argument given in Figure~\ref{fig:classes}. However, it is clear that $r\in H^-$ and yet $\alpha$ is not in $I_r$ since indeed $\alpha$ is also an equivariant Euler class associated with a wall close to $r$. This example demonstrates the need for the extended flows introduced in Section~\ref{se:equivariantMorsetheory}.

\end{document}